\documentclass[12pt,reqno]{amsart}
\usepackage{amsthm}
\usepackage{amsfonts}

\newtheorem{theorem}{Theorem}

\newtheorem{lemma}{Lemma}

\newtheorem{conj}{Conjecture}

\theoremstyle{definition}

\theoremstyle{plain}
\def\reals{{\mathbb R}}
\def\rationals{{\mathbb Q}}
\def\integers{{\mathbb Z}}
\def\naturals{{\mathbb N}}

\def\supnorm#1#2{\parallel\hspace{-.1cm}#1\hspace{-.08cm}\parallel_{#2}}
\def\Res#1{\mbox{Res}_{#1}}

\begin{document}
\title{Consequences of the Continuity of the Monic Integer Transfinite Diameter}
\author{Jan Hilmar}



\begin{abstract}
We consider the problem of determining the monic integer
transfinite diameter
$t_M(I)$ 
for real intervals $I$ of length less than $4$. We show that
$t_M([0,x])$, as a function in $x>0$, is continuous, therefore
disproving two conjectures due to Hare and Smyth. Consequently,
for $n>2\in\naturals$, we define the quantity
\begin{eqnarray*}
b_{\max}(n)&=&\sup_{b>\frac{1}{n}}\left\{b\left|t_M([0,b])=\tfrac{1}{n}\right.\right\}
\end{eqnarray*}
and give lower and upper bounds of $b_{\max}(n)$. Finally, we  improve the lower bound for
$b_{\max}(n)$ for $3\leq n\leq 8$.
\end{abstract}

\maketitle

\pagestyle{myheadings} \markboth{\textit{Monic integer transfinite
diameter}}{\textit{Jan Hilmar}}

\section{Introduction}

Let $I\subset\reals$ be a closed interval of length less than $4$ and ${\mathcal M}_n[x]$ be the set of monic polynomials of degree $n$ with integer coefficients. We define the {\it monic integer transfinite diameter} $t_M(I)$ of $I$ to be the quantity

\begin{eqnarray}
t_M(I)&=&\lim_{n\rightarrow\infty}\inf_{p_n\in{\mathcal M_n}}\supnorm{p_n}{I}^{1/n}.\label{tM}
\end{eqnarray}

Here $\supnorm{p_n}{I}=\sup_{x\in I}|p_n(x)|$ is the supremum norm of the polynomial $p_n(x)$. The problem of determining the monic integer transfinite diameter was first tackled by Borwein, Pinner and Pritsker in \cite{BPP03}. Their techniques were further developed by Hare and Smyth in \cite{SmHa05}. The problem is intimately connected to the problem of determining $t_{\integers}(I)$, the {\it integer transfinite diameter}, where the condition that the polynomials be monic is removed. Interestingly, removing this condition makes the problem much harder, as no exact values of $t_{\integers}(I)$ are known, but $t_M(I)$ can be computed explicitly in some cases.
The following lemma is an essential tool in doing so.

\begin{lemma}[\cite{BPP03}]
Let $q(x)=a_0+\cdots+a_dx^d\in\integers_n[x]$ be an irreducible polynomial with $a_d>1$ and all roots in the closed interval $I\subset\reals$ of finite length. Further, assume that $p_n(x)\in{\mathcal M}_n[x]$. Then
\begin{eqnarray*}
a_d^{-\frac{1}{d}}&\leq&\supnorm{p_n}{I}^{\frac{1}{n}}.
\end{eqnarray*}
\end{lemma}

The proof of this essentially classical result can be found in \cite{BPP03} or \cite{SmHa05} and will be omitted here.

As a consequence, $a_d^{-1/d}\leq t_M(I)$, so that polynomials
$q(x)$ as in the lemma are used to determine lower bounds on
$t_M(I)$. As a consequence, they are called {\it obstruction
polynomials} for $t_M(I)$ with {\it obstruction}  $a_d^{-1/d}$.
Since all obstructions give a lower bound, it is of interest to
find the supremum
\begin{eqnarray*}
m(I)&=&\sup\left\{\left.a_d^{-1/d}\right| q(x)=a_dx^d+\cdots +a_0, a_d>1\right\}.
\end{eqnarray*}
Here the supremum is taken over all polynomials with integer coefficients and all roots in the interval $I$. If the supremum is attained, $m(I)$ is called the {\it maximal obstruction} for $I$.

Suppose now we have an interval $I$ with maximal obstruction
$m(I)$ and find $p_n(x)\in{\mathcal M_n}$ with
$\supnorm{p_n}{I}=m(I)^n$. In this case, $m(I)\geq t_M(I)\geq
m(I)$, so that we have determined an exact value for $t_M(I)$.
Such $p_n(x)$, if it exists, is said to {\it attain the maximal
obstruction}. Some examples of this situation are as follows:

\begin{enumerate}
\item If $I=[0,1]$, then $\frac{1}{2}$ is the maximal obstruction
by $q(x)=2x-1$. At the same time,
$\supnorm{x(1-x)}{I}=\frac{1}{4}$, so that $t_M(I)=\frac{1}{2}$.
\item For an integer $n>1$, consider $I_n=[0,\frac{1}{n}]$. Then
$\frac{1}{n}$ is the maximal obstruction by $q(x)=nx-1$. At the
same time, $\supnorm{x}{I_n}=\frac{1}{n}$, so that
$t_M(I_n)=\frac{1}{n}$.
\end{enumerate}

These are just some examples to illustrate the technique. A more complete list of known values of $t_M(I)$ can be found in \cite{BPP03} and \cite{SmHa05}.

It was shown in \cite{SmHa05} that the maximal obstruction is not always attained by some $p_n(x)$ and explicit conditions for when it cannot be attained were given. The authors conjecture, however, that $t_M(I)=m(I)$ for all $I$. That this is not the case is a consequence of the continuity of a particular function, proved in Section \ref{continuity}.

\section{Continuity of $t_M(x)$, $x\geq 0$}\label{continuity}

In \cite{SmHa05}, the authors consider intervals of the form
$I=[0,b]$, where, for $1<n\in\naturals$, we have
$\frac{1}{n}<b<\frac{1}{n-1}$. From $q(x)=nx-1$, we know that
$t_M([0,b])\geq\frac{1}{n}$ and equality holds in a neighbourhood
to the right of $\frac{1}{n}$ (see Theorem \ref{bmaxbounds}).
Much more interesting is the behaviour of the function
\begin{eqnarray}
t_M(x)&=&t_M([0,x]) \label{tMx}
\end{eqnarray}
for $x\ge 0$ to the left of $\frac{1}{n}$, $n>1$. Hare and Smyth suspected that the function  had a discontinuity at $x=\frac{1}{n}$, which would agree with their conjecture that $m(I)=t_M(I)$.

To study the behaviour of $t_M(x)$, it is useful to look back at the classical paper \cite{BoEr96} of Borwein and Erd{\'e}lyi in the theory of the (non-monic) transfinite diameter. In this paper, the authors define the function $t_{\integers}(x)$ in the equivalent way and state that this function is continuous, though without the details of the proof. We will now provide the details for $t_M(x)$.

Let $T_n(x)$ be the $n^{th}$ Chebyshev polynomial on $[-1,1]$,
defined by
\begin{eqnarray}
T_n(x)&=&\cos(n\arccos x)\label{Tn}\,.
\end{eqnarray}
This can be rewritten as
\begin{eqnarray*}
T_n(x)&=&\tfrac{1}{2}\left[\left(x+\sqrt{x^2-1}\right)^n+\left(x-\sqrt{x^2-1}\right)^n\right].
\end{eqnarray*}
From this it immediately follows that
\begin{eqnarray}
T_n(x)&\leq&\left(x+\sqrt{x^2-1}\right)^n\mbox{ for }x\geq
1\,.\label{Tnineq}
\end{eqnarray}

We will also need Chebyshev's inequality from \cite{BorErd}:

\begin{lemma}\label{Cheb}
Let $q\in\reals[x]$. Then, for $x\in\reals\backslash[-1,1]$,
\begin{eqnarray}
\bigl|q(x)\bigr|&\leq&\bigl|T_n(x)\bigr|\supnorm{q}{[-1,1]}\,.
\end{eqnarray}
\end{lemma}

We can then prove:

\begin{lemma}\label{lemmanormineq}
Let $b>b_0>0$, $p_n\in\reals_n[y]$. Then, for every $\delta>0$, there exists $k_{b,\delta}$, not depending on $n$, such that
\begin{eqnarray}
\supnorm{p_n}{[0,b+\delta]}&\leq&(1+k_{b,\delta})^n\supnorm{p_n}{[0,b]}\,,\label{normineq}
\end{eqnarray}
with $\lim_{\delta\rightarrow0}k_{b,\delta}=0$ for fixed $b$.
\end{lemma}

\begin{proof}
Given $p_n\in\reals_n[y]$, let $y\in[0,b]$ and $x=\frac{2}{b}y-1$. Then $x\in[-1,1]$. Put $q_n(x)=p_n(y)$. Then, by Lemma \ref{Cheb}, for $x\not\in [-1,1], y\not\in [0,b]$, we have
\begin{eqnarray*}
|p_n(y)|=|q_n(x)|&\leq&\left|T_n(x)\right|\supnorm{q_n}{[-1,1]}\\
&=&\left|T_n\left(\tfrac{2}{b}y-1\right)\right|\supnorm{p_n}{[0,b]}.
\end{eqnarray*}
Note also that 
\begin{eqnarray*}
\max_{y\in[b,b+\delta]}\left|T_n\left(\tfrac{2}{b}y-1\right)   \right|&=&\max_{x\in[1,1+2\frac{\delta}{b}]}\left|T_n(x)\right|\\
&=&\supnorm{T_n}{[1,1+2\frac{\delta}{b}]}\,.
\end{eqnarray*}
This clearly implies that
\begin{eqnarray*}
\supnorm{p_n}{[b,b+\delta]}&\leq&\supnorm{T_n}{[1,1+2\frac{\delta}{b}]}\supnorm{p_n}{[0,b]}.
\end{eqnarray*}
Using inequality (\ref{Tnineq}) above, we see that
\begin{eqnarray*}\supnorm{T_n}{[1,1+2\frac{\delta}{b}]}&\leq&\left(1+2\tfrac{\delta}{b}\left(1+\sqrt{1+\tfrac{b}{\delta}}\right)\right)^n\,.
\end{eqnarray*}

The result now follows by letting $k_{b,\delta}=2\frac{\delta}{b}\left(1+\sqrt{1+\frac{b}{\delta}}\right)>0$ and observing that
\begin{eqnarray*}
\supnorm{p_n}{[0,b+\delta]}&=&\max\left\{\supnorm{p_n}{[0,b]},\supnorm{p_n}{[b,b+\delta]}\right\}\\
&\leq&\max\left\{\supnorm{p_n}{[0,b]},(1+k_{b,\delta})^n\supnorm{p_n}{[0,b]}\right\}\\
&=&(1+k_{b,\delta})^n\supnorm{p_n}{[0,b]}.
\end{eqnarray*}
\end{proof}

Using this inequality, we also get that, for $b,\delta>0$ fixed,
\begin{eqnarray}
\supnorm{p_n}{[0,b-\delta]}&\geq&\supnorm{p_n}{[0,b]}\left(\frac{1}{1+k_{b-\delta,\delta}}\right)^n.\label{invnormineq}
\end{eqnarray}
Note also that $\lim_{\delta\rightarrow0}k_{b-\delta,\delta}=0$.

We can now use this to prove
\begin{theorem}The function \label{tMcontinuity}
$t_M(x)$ is continuous on $(0,\infty)$.
\end{theorem}

\begin{proof}
First, note that $t_M(x)$ is (non-strictly) increasing in $x$. Let $b\in(0,\infty)$, $\epsilon>0$ and choose $\delta=\min\{\delta_1,\delta_2\}$, where $\delta_1$ is chosen such that $k_{b,\delta_1}<\frac{\epsilon}{t_M(b)}$ and $\delta_2$ is such that $\frac{k_{b-\delta_2,\delta_2}}{1+k_{b-\delta_2,\delta_2}}<\frac{\epsilon}{t_M(b)}$.

Let $0<|b-x|<\delta$. The argument splits into two cases:

\noindent (1) Suppose that $0<b-x<\delta\leq\delta_1$. Since $t_M(x)$ is
increasing, we have
\begin{eqnarray*}
0\leq t_M(x)-t_M(b)&\leq& t_M(b+\delta_1)-t_M(b)\\
&=&\lim_{n\rightarrow\infty}\left(\inf_{p_n\in{\mathcal M}_n[x]}\supnorm{p_n}{[0,b+\delta_1]}^{1/n}-\inf_{p_n\in{\mathcal M}_n[x]}\supnorm{p_n}{[0,b]}^{1/n}\right)\\
&\leq&\lim_{n\rightarrow\infty}\left(\inf_{p_n\in{\mathcal M}_n[x]}k_{b,\delta_1}\supnorm{p_n}{[0,b]}^{1/n}\right)\\
&=&t_M(b)k_{b,\delta_1}<\epsilon\,.
\end{eqnarray*}

\noindent (2) Now assume that $0<x-b<\delta\leq\delta_2$. Here, we
get
\begin{eqnarray*}
0\leq t_M(b)-t_M(x)&\leq& t_M(b)-t_M(b-\delta_2)\\
&=&\lim_{n\rightarrow\infty}\left(\inf_{p_n\in{\mathcal M}_n[x]}\supnorm{p_n}{[0,b]}^{1/n}-\inf_{p_n\in{\mathcal M}_n[x]}\supnorm{p_n}{[0,b-\delta_2]}^{1/n}\right)\\
&\leq&\lim_{n\rightarrow\infty}\left(\inf_{p_n\in{\mathcal M}_n[x]}\left(\frac{k_{b-\delta_2,\delta_2}}{1+k_{b-\delta_2,\delta_2}}\right)\supnorm{p_n}{[0,b]}^{1/n}\right)\\
&=&t_M(b)\frac{k_{b-\delta_2,\delta_2}}{1+k_{b-\delta_2,\delta_2}}<\epsilon\,.
\end{eqnarray*}
Thus, for $0<|b-x|<\delta$, we have $|t_M(b)-t_M(x)|<\epsilon$ for any $b\in(0,\infty)$, proving continuity for $x>0$.
\end{proof}

As mentioned before, Borwein and Erd{\'e}lyi stated this result
for the (non-monic) integer transfinite diameter. In fact, if ${\mathcal
A_n}[x]\subseteq\reals_n[x]$ and
\begin{eqnarray}
t_{{\mathcal A}}(I)&=&\lim_{n\rightarrow\infty}\inf_{0\not\equiv
p_n\in{\mathcal A_n[x]}}\supnorm{p_n}{I}^{\frac{1}{n}}\,,
\end{eqnarray}
one can define $t_{{\mathcal A}}(x)$ in the equivalent way and the prove continuity of this function for $x\geq 0$ as in Theorem \ref{tMcontinuity}.

The continuity of $t_M(x)$ sheds some light on conjectures made by Hare and Smyth in \cite{SmHa05}:

\begin{conj}[Zero-endpoint interval conjecture]\label{zepconj}
If $I=[0,b]$ with $b\leq 1$, then $t_M(I)=\frac{1}{n}$, where $n=\max\left(2,\lceil\frac{1}{b}\rceil\right)$ is the smallest integer $n\geq 2$ for which $\frac{1}{n}\leq b$.
\end{conj}

\begin{conj}[Maximal obstruction implies $t_M(I)$ conjecture]\label{mitconj}
If an interval $I$ of length less than $4$ has a maximal obstruction $m(I)$, then $t_M(I)=m(I)$.
\end{conj}

Theorem \ref{tMcontinuity} clearly shows Conjecture \ref{zepconj} to be false: since $t_M(x)$ is continuous to the left of $\frac{1}{n}, n\in\naturals$ and $t_M([0,\frac{1}{n}])=\frac{1}{n}$, we cannot have $t_M([0,b])=\frac{1}{n+1}$ for all $b<\frac{1}{n}$, as claimed in the conjecture. A further implication is that for these intervals, $t_M(I)\neq m(I)$, contrary to Conjecture \ref{mitconj}.

\section{The function $b_{\max}(n)$}
It turns out that $t_M(x)$ is indeed constant on a large interval
to the right of $\frac{1}{n+1}, n\in\naturals$. We define for
$n>1\in\naturals$,
\begin{eqnarray}
b_{\max}(n)&=&\sup_{b>\frac{1}{n}}\left\{b\left| t_M(b)=\tfrac{1}{n}\right.\right\}.\label{bmax}
\end{eqnarray}

For $n=1$, this quantity is not finite, as $t_M(I)=1$ for $|I|\geq 4$ (see \cite{BPP03} for details). For $n=2$, we can use the results in \cite{SmHa05} to obtain $1.26\leq b_{\max}(2)<1.328$. For $n>2$, we have the following:

\begin{theorem}\label{bmaxbounds}
Let $n>2\in\naturals$. Then
\begin{eqnarray*}
\frac{1}{n}+\frac{1}{n^2(n-1)}&<
b_{\max}(n)\leq&\frac{4n}{(2n-1)^2}\,.
\end{eqnarray*}
\end{theorem}

\begin{proof}
The first inequality follows from the polynomial
\begin{eqnarray*}
P_n(x)&=&x^{n^2-2}(x^2-nx+1)\,.
\end{eqnarray*}

This polynomial,  first used in \cite{SmHa05}, was shown to
have the following properties:
\begin{enumerate}
\item $P_n(\frac{1}{n})=\left(\frac{1}{n}\right)^{n^2}$.
 \item
$P_n'(\frac{1}{n})=0$ and the polynomial  has no other extrema in
$[0,\frac{1}{n-1}]$. \item $P_n(x)$ has a root
$\beta_n=\frac{2}{n+\sqrt{n^2-4}}>\frac{1}{n}$, and $|P_n(x)|$ is
strictly increasing in $(\beta_n,\frac{1}{n-1})$.
\end{enumerate}
These properties were used in \cite{SmHa05} to show that $\supnorm{P_n}{[0,\frac{1}{n}+\epsilon]}=\left(\frac{1}{n}\right)^{n^2}$ for some $\epsilon>0$.

Evaluating $P_n(x)$ at $x=\frac{1}{n}+\frac{1}{n^2(n-1)}$ gives
\begin{eqnarray*}
\left|P_n\left(\frac{1}{n}+\frac{1}{n^2(n-1)}\right)\right|&=&\left(\frac{n^2-n+1}{n^2-n}\right)^{n^2}\frac{n^3-3n^2+2n-1}{(n^2-n+1)^2}\,.\\
\end{eqnarray*}
To show that this is indeed less than $(\frac{1}{n})^{n^2}$, first
note that the sequence
\begin{eqnarray*}
\left\{\left(\frac{n^2-n}{n^2-n+1}\right)^{n^2}\right\}_{n=1}^{\infty}
\end{eqnarray*}
is increasing and tends to $e^{-1}$. As it is increasing, we
clearly have
\begin{eqnarray}
\left(\frac{n^2-n}{n^2-n+1}\right)^{n^2}&\geq&\left(\frac{2}{3}\right)^4\mbox{
for }n>2\,.\label{seqineq}
\end{eqnarray}

Further, note that for all $n$,
\begin{eqnarray}
\left(\frac{2}{3}\right)^4&>&\frac{n^3-3n^2+2n-1}{(n^2-n+1)^2}\,.\label{ratfunineq}
\end{eqnarray}
Thus taking (\ref{seqineq}) and (\ref{ratfunineq}) together, we
have, for $n>2$,
\begin{eqnarray*}
\left(\frac{n^2-n}{n^2-n+1}\right)^{n^2}&>&\frac{n^3-3n^2+2n-1}{(n^2-n+1)^2}\,.
\end{eqnarray*}
Rearranging now gives the desired result.

For the upper bound, one has to look directly at (\ref{normineq}).
Suppose we have some $p_d(x)\in{\mathcal M_d}[x]$ such that
$\supnorm{p_d}{I_{\delta_n}}^{\frac{1}{d}}=\frac{1}{n}$ on an
interval $I_{\delta_n}=[0,\frac{1}{n-1}-\delta_n]$. Clearly,
$\supnorm{p_d}{[0,\frac{1}{n-1}]}^{\frac{1}{d}}\geq\frac{1}{n-1}$
since $\frac{1}{n-1}\leq t_M([0,\frac{1}{n-1}])$. Thus, using
(\ref{normineq}), we get
\begin{eqnarray*}
\frac{1}{n-1}&\leq&\frac{1}{n}\bigl(1+k_{\frac{1}{n-1}-\delta_n,\delta_n}\bigr)\,.
\end{eqnarray*}

Using the explicit expression for
$k_{\frac{1}{n-1}-\delta_n,\delta_n}$ obtained in the proof of
Lemma~\ref{lemmanormineq}, we see that then
$\delta_n\geq\delta_{\min}(n)=\frac{1}{4n^3-8n^2+5n-1}$, thus
obtaining
\begin{eqnarray*}
b_{\max}(n)&\leq&\frac{1}{n-1}-\delta_{\min}(n)=\frac{4n}{(2n-1)^2}\,.
\end{eqnarray*}
\end{proof}

Using the computational methods outlined in Section
\ref{computational}, we get improved lower bounds for $b_{\max}(n)$ for $n=3,\ldots,8$. This is done by finding a $b\in(\frac{1}{n},\frac{1}{n-1})$ as large as possible and a polynomial $P_n(x)$ with $\supnorm{P_n}{[0,b)}=\frac{1}{n}$, so that then $b_{\max}(n)\ge b.$
 The polynomials $P_n$ are given in  Table \ref{bmaxpolys}. The polynomial $P_3$ is a corrected version of  one appearing in \cite{SmHa05}, which does not have the property claimed, while $P_4$ appears in \cite{SmHa05}.
\begin{eqnarray*}
\begin{array}{cc}
0.465\leq b_{\max}(3)\,,&0.303\leq b_{\max}(4)\,,\\
0.230\leq b_{\max}(5)\,,&0.184\leq b_{\max}(6)\,,\\
0.148\leq b_{\max}(7)\,,&0.130\leq b_{\max}(8)\,.
\end{array}
\end{eqnarray*}

As $n$ gets larger, computations become increasingly difficult, as the difference $\frac{1}{n-1}-\frac{1}{n}$ becomes small compared to $\frac{1}{n}$.

Using (\ref{normineq}), one can obtain a new lower bound for $t_M([0,b]), b<1$:

\begin{lemma}\label{extralb}
Let $I_b=[0,b]$, $b<1$ and let $n=\min\{m\in\naturals\mid \frac{1}{m}>b\}$. Then
\begin{eqnarray*}
t_M(I_b)&\geq&\max\left\{\frac{1}{n+1},\frac{b}{2(1+\sqrt{1-nb})-nb}\right\}\,.
\end{eqnarray*}
\end{lemma}

\begin{proof}
Let $\delta=\frac{1}{n}-b$. As can easily be seen from
(\ref{normineq}),
\begin{eqnarray*}
t_M\left(\left[0,\tfrac{1}{n}-\delta\right]\right)&\geq&t_M\left(\left[0,\tfrac{1}{n}\right]\right)\frac{1}{1+k_{\frac{1}{n}-\delta,\delta}}\\
&=&\frac{1-n\delta}{n(1+n\delta+2\sqrt{n\delta})}\\
&=&\frac{b}{2-nb+2\sqrt{1-nb}}\,.
\end{eqnarray*}
Seeing that $\frac{1}{n+1}$ is a larger lower bound for $b\leq
b_{\text{max}}(n+1)$, we get the result.
\end{proof}

\section{The Farey interval conjecture}
Another open conjecture, this one taken from \cite{BPP03}, is the following:

\begin{conj}[Farey Interval Conjecture]
Let $\frac{p}{q},\frac{r}{s}\in\rationals$ with $q,s>0$ be such
that $rq-ps=1$. Then
\begin{eqnarray*}
t_M\left(\left[\frac{p}{q},\frac{r}{s}\right]\right)&=&\max\left\{\frac{1}{q},\frac{1}{s}\right\}.
\end{eqnarray*}
\end{conj}

Computationally, the authors verified the conjecture for
denominators up to $21$ and it was proved for an infinite family
of such intervals in \cite{SmHa05}. It is perhaps worth noting
that continuity of $t_M(x)$ cannot be used to find a
counterexample to this conjecture, as the following argument
shows.

 Let $n>1\in\naturals$. We will show that we cannot find a
Farey Interval of the form $[\frac{k}{n},\frac{p}{q}]$ with $1\leq
k<n, n=\min\{q,n\}$ and $\frac{p}{q}>b_{\max}^*(n)$, where
\begin{eqnarray*}
b_{\max}^*(n)&=&\sup_{\frac{k}{n}<b}\left\{b\left|
t_M\left(\left[\tfrac{k}{n},b\right]\right)=\tfrac{1}{n}\right.\right\}\,.
\end{eqnarray*}

As can easily be derived from the proof of Lemma
\ref{lemmanormineq},
\begin{eqnarray*}
b^*_{\max}(n)&\leq&\frac{k(4n^2+1)}{n(2n-1)^2}\,.
\end{eqnarray*}
If we wanted to use this to derive a counterexample to the Farey
Interval Conjecture, we would need $\frac{p}{q}>b^*_{\max}(n)$.
Using the Farey property
\begin{eqnarray}
pn-qk=1\,,\label{Fareyprop}
\end{eqnarray}
we can write this as
\begin{eqnarray*}
\frac{1+qk}{qk}&>&\frac{4n^2+1}{(2n-1)^2}\,.
\end{eqnarray*}

From this it follows that $1+qk<2n+1$, so that, using (\ref{Fareyprop}) again, $1+qk=2n$. Now, as $q>n$, it is clear that $k=1$ for this to hold.

In the case $k=1$, one can show that the Farey interval is then of
the form $\bigl[\frac{1}{n},\frac{1+t}{(1+t)n-1}\bigr]$,
$t\in\naturals$. But no such interval with
$\frac{1+t}{(1+t)n-1}>b^*_{\max}(n)=\frac{4n^2+1}{(2n-1)^2n}$
exists.

The result for the remaining Farey intervals is obtained by using the transformations $x\mapsto m\pm x, m\in\integers$.

\section{Computational Methods}\label{computational}

In order to improve the lower bounds for $b_{\max}(n)$ given in Theorem \ref{bmaxbounds}, we need to turn to computational methods to attempt to find a monic polynomial $P(x)\in\integers[x]$ attaining the maximal obstruction on an interval $[0,b)$ with $\frac{1}{n}<b<\frac{1}{n-1}$. 
These come in two stages:

\begin{enumerate}
\item Using a modification of the LLL algorithm to find factors $f_i(x)$ of $P(x)$.
\item Using Linear Programming methods first used in \cite{BoEr96} in connection with the integer transfinite diameter with additional equality constraints obtained in \cite{SmHa05} to determine the exponents $\alpha_i$.
\end{enumerate}

We will briefly discuss the implementations of both parts of the algorithm.

\begin{enumerate}
\item LLL is an algorithm that, given a basis ${\mathbf b}$ for a
lattice $\Lambda$, produces a `small' basis for $\Lambda$ with
respect to a given inner product $\left<\cdot,\cdot\right>$. In
their modification of the LLL algorithm for monic polynomials
introduced in \cite{BPP03}, the authors used the Lattice
$\integers_n[x]$ with the basis ${\mathbf b}=(1,x,x^2,\cdots x^n)$
and the inner product
\begin{eqnarray*}
\left<p_n,q_n\right>&=&\int_a^bp_n(x)q_n(x)dx+a_nb_n
\end{eqnarray*}
for $p_n(x)=a_nx^n+\cdots+a_0, q_n(x)=b_nx^n+\cdots+b_0\in\integers_n[x]$. The additional factor $a_nb_n$ is used to discourage non-monic factors from appearing, and the algorithm usually produces only one monic basis element of degree $n$.

In practice, we used the following recursive algorithm to identify factors $f_i(x)$ of $P(x)$ for an interval $I=[a,b]$ where the maximal obstruction polynomial $q(x)=a_dx^d+\cdots+a_0$ is known:

\begin{enumerate}
\item Start with ${\mathbf b}=(1,x,x^2,\cdots,x^k)$ for $k=20$ (in
some cases, a larger basis was required initially). \item Run LLL,
generating a list of factors $l=\{f_i(x)\}$. \item Sieve the list
by using the condition that if $f_i(x)\mid P(x)$, then the
resultant has to satisfy $|\Res{}(f_i,q)|=1$ (see \cite{SmHa05}).
\item For every $f_i$ still in $l$, define $$\widehat{{\mathbf
b}}_i=(1,f_i(x),f_i(x)x,f_i(x)x^2,\cdots f_i(x)x^k)$$ and re-run
the LLL Algorithm with this basis, adding new factors to $l$.
\item Repeat steps (a)--(d) until no more new factors are found,
at which point we return $l$.
\end{enumerate}

\item To determine the exponents $\alpha_i$ of $f_i(x), 1\leq
i\leq N$, we use a technique first used by Borwein and Erd{\'e}lyi
in \cite{BoEr96}. Given a list of factors $l=\{f_i(x)\}$, one
attempts to minimise $m$ subject to
\begin{eqnarray}
\begin{array}{l}
\text{(i)}\, \sum_{i=1}^N\frac{\alpha_i}{\deg f_i}\log
|f_i(x)|\leq m-g(x)\,,\quad\text{for}\,
x\in X\,,\\ \\
\text{(ii)}\,\sum_{i=1}^N\alpha_i=1\,,\\ \\
\text{(iii)}\,\sum_{i=1}^N\frac{\alpha_i}{\deg
f_i}\frac{f'(\beta_s)}{f(\beta_s)}=0\,,\quad\text{for}\, 1\leq
s\leq \deg q\,,\text{where}\,
q(\beta_s)=0\,,\\ \\
\text{(iv)}\, \alpha_i\geq 0\qquad\text{for}\, 1\leq i\leq N\,,
\end{array}\label{consts}
\end{eqnarray}
over a finite set $X\subset I$. Here, $g(x)$ is a function such
that
\begin{eqnarray*}
g(x)&=&\left\{
\begin{array}{ll}
0\,,&q(x)=0\,,\\
\epsilon(x)>0\,,&q(x)\neq 0\,.
\end{array}\right.
\end{eqnarray*}

The use of this function is theoretically not necessary, but is
useful when doing computations, as it avoids having to deal with
exact values at points where the polynomial does not need to
attain the maximal obstruction.

The first two constraints in (\ref{consts}) are taken from \cite{BoEr96} with a slight modification to the first, while the third is unique to the monic case and taken from \cite{SmHa05}. This is also where we get the final set of constraints:

Let $\beta_s$ be a root of $q(x)$ and define
$\hat{f}^{(s)}_i=\frac{1}{\deg f_i}\log|f_i(\beta_s)|$. If
$b_1=-\frac{1}{d}\log |a_d|,b_2,\ldots, b_l$ is an independent
generating set for the $\integers$-lattice generated by
$-\frac{1}{d}\log|a_d|$ and the $\hat{f}^{(s)}_i$, let
$c^{(s)}_{j,i}$ be such that
\begin{eqnarray*}
\sum_{j=1}^lc^{(s)}_{j,i}b_j&=&\hat{f}^{(s)}_i.
\end{eqnarray*}
Then we get the additional conditions, derived in \cite{SmHa05}:
\begin{eqnarray}
\sum_{i=1}^Nc^{(s)}_{j,i}\alpha_i&=&\left\{
\begin{array}{rl}
-\frac{1}{\deg q}\,,&j=1\,,\\
0\,,&j>1\,,
\end{array}\right.\qquad\text{ for } 1\leq s\leq\deg q\label{Harecond}
\end{eqnarray}

Again, we use a recursive algorithm for determining the exponents.
Given a set of points $X_k$, we use (\ref{consts}) and
(\ref{Harecond}) to determine the optimal exponents
$\{\alpha_1^{(k)},\alpha_2^{(k)},\ldots,\alpha_N^{(k)}\}$
attaining the minimum value $m_k$. Then, we construct the
normalised `polynomial'
\begin{eqnarray*}
P^{(k)}(x)&=&\prod_{i=1}^Nf_i(x)^{\frac{\alpha_i^{(k)}}{\deg
f_i}}\,,
\end{eqnarray*}
and add its extrema to $X_k$ to obtain $X_{k+1}$. Starting with a small set of values $X_1\subset I$, we repeat this procedure until we get $K\in\naturals$ such that $|m_K-m_{K-1}|<\epsilon$ for required precision $\epsilon>0$.

Finally, we compute the supremum norm of $P_K(x)\approx e^{m_K}$ on the interval and verify that $\supnorm{P_K}{I}=|a_d|^{-\frac{1}{d}}$. One can attempt to find rational approximations of smaller denominator to the exponents, always checking that the obstruction is still attained. The attaining polynomial $P(x)$ is then found by clearing denominators in the exponents of $P_K(x)$.
\end{enumerate}

\noindent\textbf{Acknowledgement}

 The author would like to thank
Kevin Hare for helpful insights and Chris Smyth for guidance and
relentless proofreading. Further, I would like to thank the
referee for carefully reading the manuscript and catching some
mistakes as well as providing useful comments. \vfill\eject

\pagebreak
\begin{table}
\caption{\label{bmaxpolys}Polynomials used for lower bounds on $b_{\max}(n)$}
{\footnotesize
\begin{eqnarray*}
\begin{array}{|  l |}
\hline
\\
\begin{array}{l}P_3(x)=x^{45944640}(x^{14}-11406261x^{13}+47054086x^{12}-88456310x^{11}\\+100247244x^{10}-76341256x^9+
 41208853x^8-16202606x^7+4692047x^6-999261x^5\\+154318x^4-16766x^3+1211x^2
  -52x+1)^{2450525}\\
(x^8+14184x^7-34944x^6+36442x^5-20832x^4+7041x^3-1405x^2+153x-7)^{877415}\\
(x^8+4842x^7-10935x^6+10355x^5-5317x^4+1594x^3-278x^2+26x-1)^{2571030}\\
(x^8+7812x^7-18072x^6+17561x^5-9271x^4+2864x^3-516x^2+50x-2)^{595980}\\
(x^7-1233x^6+2406x^5-1913x^4+791x^3-179x^2+21x-1)^{1210840}\\
(x^5-3x^4+7x^3-11x^2+6x-1)^{1052898}\\
\end{array}\\
\\
\hline
\\
\begin{array}{l}
P_4(x)=x^{640}(x^5+432x^4-456x^3+179x^2-31x+2)^{47}\\
(x^7+8760x^6-13342x^5+8488x^4-2784x^3+514x^2-50x+2)^{35}
\end{array}\\
\\
\hline
\\
\begin{array}{l}
P_5(x)=x^{1050990}(x^{10}+5544095x^9-9115714x^8+6623719x^7-2790988x^6\\+751349x^5-133974x^4+15818x^3-1192x^2+52x-1)^{78796}\\(x^6+4950x^5-4605x^4+1698x^3-310x^2+28x-1)^{21825}
\end{array}\\
\\
\hline
\\
\begin{array}{l}
P_6(x)=x^{5232473}(x^5+1260x^4-852x^3+215x^2-24x+1)^{118824}\\(x^7-140190x^6+132517x^5-51966x^4+10819x^3-1261x^2+78x-2)^{200917}
\end{array}\\
\\
\hline
\\
\begin{array}{l}P_7(x)=x^{44}(x^5+3472x^4-1826x^3+358x^2-31x+1)\end{array}\\
\\
\hline
\\
\begin{array}{l}P_8(x)=x^{12288}(x^2-8x+1)^{246}(x^4-576x^3+208x^2-25x+1)^{741}\end{array}\\
\\
\hline
\end{array}
\end{eqnarray*}}
\end{table}


\end{document}